\documentclass[10pt]{amsart}
\usepackage{amssymb, amscd, amsmath, amsthm, amscd,
epsf, epsfig, latexsym, enumerate, pstricks,graphics}

\begin{document}

\title{an explicit formula for a branched covering from $\mathbb{CP}^2$ to $S^4$}
 
\author{J.A.Hillman}
\address{School of Mathematics and Statistics\\
     University of Sydney, NSW 2006\\
      Australia }

\email{jonathan.hillman@sydney.edu.au}

\begin{abstract}
We give an explicit formula for a 2-fold branched covering from $\mathbb{CP}^2$ to $S^4$,
and relate it to other maps between quotients of $S^2\times{S^2}$. 
\end{abstract}



\maketitle

It is well known that the quotient of the complex projective plane 
$\mathbb{CP}^2$ by complex conjugation is the 4-sphere \cite{Ku,Ma}.
(See also \cite{AK}.)
A key element of Massey's exposition is the identification of 
$\mathbb{CP}^2$ with the 2-fold symmetric product $S^2(\mathbb{CP}^1)$, 
the quotient of $\mathbb{CP}^1\times\mathbb{CP}^1$ by
the involution which exchanges the factors.
Lawson has displayed very clearly the topology underlying these facts \cite{La}.
We shall use harmonic coordinates for real and complex projective spaces
to give explicit formulae for some of these quotient maps.
(We describe briefly the work of Massey and Kuiper at the end of this note.)

A smooth map $f:M\to{N}$ between closed $n$-manifolds is a 
2-fold branched covering if $M$ has a codimension-2 submanifold  $B$
(the branch locus), such that $f|_{M\setminus{B}}$ is a 2-to-1 immersion,
$f|_B:B\to{f(B)}$ is a bijection onto a submanifold $f(B)$ (the branch set)
and along $B$ the map $f$ looks like $(b,z)\mapsto(f(b),z^2)$ in local coordinates, with $b\in{B}$ and 
transverse complex coordinate $z\in\mathbb{C}$.

We shall view $S^2$ as the unit sphere in $\mathbb{C}\times\mathbb{R}$.
The antipodal involution $A$ is given by multiplying the coordinates by $-1$.
Let $\Delta_{S^2}$ and $\Gamma_A$ be the diagonal of $S^2\times {S^2}$ and the graph of $A$,
respectively.
Let $\sigma$ and $\tau$ be the diffeomorphisms of $S^2\times{S^2}$
given by $\sigma(s,s')=(s',A(s))$ and $\tau(s,s')=(s',s)$,
for $s,s'\in{S^2}$.
Then $\sigma$ and $\tau$ generate a dihedral group of order 8,
since $\sigma^4=\tau^2=1$ and $\tau\sigma\tau=\sigma^{-1}$.
This group acts freely on the complement of $\Delta_{S^2}\cup\Gamma_A$.
In particular, $\sigma$ acts freely on $S^2\times{S^2}$ and $\sigma(\Delta_{S^2})=\Gamma_A$,
while $\tau$ fixes $\Delta_{S^2}$ pointwise and acts freely on $\Gamma_A$.

Let $\gamma:S^2\to\mathbb{CP}^1$ be the stereographic projection,
given by $\gamma(z,t)=[z:1-t]$, for $(z,t)\in{S^2}$.
Then its inverse is given by
\[
\gamma^{-1}([u:v])=(\frac{2u\bar{v}}{|u|^2+|v|^2},\frac{|u|^2-|v|^2}{|u|^2+|v|^2}),
\]
and the antipodal map is given in harmonic coordinates  by
\[
\gamma{A}\gamma^{-1}([u:v])=[-\bar{v}:\bar{u}],
\]
for $[u:v]\in\mathbb{CP}^1$.
The image of $\Delta_{S^2}$ under $\gamma\times\gamma$ is the diagonal
$\Delta\subset\mathbb{CP}^1\times\mathbb{CP}^1$.

Let $\alpha$, $\beta$ and $f:\mathbb{C}^2\to\mathbb{C}^2$
be the maps given by $\alpha(w,z)=(\frac{w-z}2,\frac{w+z}2)$, 
$\beta(w,z)=(z,z^2-4w)$ and $f(w,z)=(wz,w+z)$,
for $(w,z)\in\mathbb{C}^2$.
Then $\alpha$ and $\beta$ are biholomorphic, and $\beta{f}\alpha(w,z)=(w,z^2)$.
Therefore $f$ is a 2-fold branched covering, branched over
$\{(p,q)\in\mathbb{C}^2\mid{q^2=4p}\}$.
The extension
$\hat{f}:\mathbb{CP}^1\times\mathbb{CP}^1\to\mathbb{CP}^2$ given by
\[\hat{f}([u:v],[u':v'])=[uu':uv'+u'v:vv']\]
is a 2-fold branched covering, branched over $\hat{f}(\Delta)$.
(The restriction $\hat{f}|_\Delta$ is essentially the 2-canonical embedding of $\mathbb{CP}^1$ in $\mathbb{CP}^2$,
with image the conic $V(z_1^2-4z_0z_2)=\{[z_0:z_1:z_2]\in\mathbb{CP}^2\mid{z_1^2=4z_0z_2}\}$.)
The composite $\lambda=\hat{f}(\gamma\times\gamma)$ is essentially Lawson's map, 
giving 
\[
S^2\times{S^2}/\langle\tau\rangle\cong{S^2(\mathbb{CP}^1)}=
\mathbb{CP}^1\times\mathbb{CP}^1/\langle\gamma\tau\gamma^{-1}\rangle\cong\mathbb{CP}^2.
\]

Let $c_n:\mathbb{CP}^n\to\mathbb{CP}^n$ be complex conjugation,
for $n\geq1$.
Then $\hat{f}(c_1\times{c_1})=c_2\hat{f}$.
Lawson observed that if $\theta:\mathbb{CP}^2\to\mathbb{CP}^2$ 
is the linear automorphism given by
\[\theta([u:v:w])=[(iu+w):(1-i)v:(u+iw)]\]
then 
\[
\lambda\sigma^2=\theta^2c_2\lambda=\theta{c_2}\theta^{-1}\lambda.
\]
Hence $c_2$ is conjugate to a map covered by the {\it free\/} involution $\sigma^2$.
(Note that $\theta^2([u:v:w])=[w:-v:u]$,  $c_2\theta{c_2}=\theta^{-1}$ and $\theta^4=id_{\mathbb{CP}^2}$.)

If we identify $\mathbb{C}$ with $\mathbb{R}^2$ and use 
{\it real\/} harmonic coordinates we may instead extend $f$ to a map 
$g:\mathbb{RP}^2\times\mathbb{RP}^2\to\mathbb{RP}^4$,
given by
\[g([r:s:t],[r':s':t'])=[rr'-ss':rs'+sr':rt'+tr':st'+ts':tt'].\]
This is a 2-fold branched covering, with branch locus the diagonal,
and induces a diffeomorphism
\[
\mathbb{RP}^2\times\mathbb{RP}^2/(x,y)\sim(y,x)\cong\mathbb{RP}^4.
\]

The map $g$ has a lift $\tilde{g}:S^2\times{S^2}\to{S^4}$, 
given by
\[
\tilde{g}((r,s,t),(r',s',t'))=\nu(rr'-ss',rs'+sr',rt'+tr',st'+ts',tt'),
\]
where $\nu:\mathbb{R}^5\setminus\{O\}\to{S^4}$ is the radial normalization.
(If $(r,s,t),(r',s',t')\in{S^2}$ then the norm of
$(rr'-ss',rs'+sr',rt'+tr',st'+ts',tt')$ is $\sqrt{1+2tt'(rr'+ss')}$.)
We may also obtain $\tilde{g}$ by normalizing the map
\[
((z,t),(z',t'))\mapsto(zz',zt'+tz',tt')\in\mathbb{C}^2\times\mathbb{R}\setminus\{O\},
\quad\forall(z,t),(z',t')\in{S^2}.
\]
This map is invariant under $\sigma^2$ and $\tau$, 
and is generically 4-to-1.
Hence it factors through a map $g^+:S^2\times{S^2}/\langle\sigma^2\rangle\to{S^4}$,
and induces maps $G=\tilde{g}\lambda^{-1}:\mathbb{CP}^2\to{S^4}$ and
$h:S^2\times{S^2}/\langle\sigma\rangle\to\mathbb{RP}^4$.

The lattice of quotients of $S^2\times{S^2}$ by the subgroups of the group $\langle\sigma,\tau\rangle$ generated by
$\sigma$ and $\tau$ is a commuting diagram:

\setlength{\unitlength}{1mm}
\begin{picture}(90,55)(-6,-3)

\put(45,45){$S^2\times{S^2}$}

\put(5,30){$S^2\times\mathbb{RP}^2$}
\put(30,30){$\mathbb{RP}^2\times{S^2}$}
\put(55,30){$S^2\times{S^2}/\langle\sigma^2\rangle$}
\put(90,30){$\mathbb{CP}^2$}

\qbezier(52,43)(33,38.5)(14,34)
\qbezier(52,43)(47.5,38.5)(43,34)
\qbezier(52,43)(57,38.5)(62,34)
\qbezier(52,43)(71,38.5)(90,34)

\put(20,15){$\mathbb{RP}^2\times\mathbb{RP}^2$}
\put(45,15){$S^2\times{S^2}/\langle\sigma\rangle$}
\put(75,15){$S^4$}

\qbezier(14,29)(22,24)(30,19)
\qbezier(39,29)(34.5,24)(30,19)
\qbezier(64,29)(58.5,24)(53,19)
\qbezier(64,29)(47,24) (30,19)

\qbezier(90,29)(83.5,24)(77,19)
\qbezier(64,29)(70.5,24)(77,19)

\put(50,0){$\mathbb{RP}^4$}

\qbezier(31,14)(41.5,9)(52,4)
\qbezier(52,14)(52,9)(52,4)
\qbezier(73.6,15)(62.8,9.5)(52,4)

\put(73,40){$\lambda$}
\put(38,6){$g$}
\put(73,24){$g^+$}
\put(85.5,22){$G$}

\put(48.5,9){$h$}

\end{picture}

The five nontrivial subgroups of $\langle\sigma,\tau\rangle$ 
that do {\it not}  contain $\tau$ (namely,   
$\langle\sigma\tau\rangle$, $\langle\tau\sigma\rangle$,
$\langle\sigma^2\rangle$, $\langle\sigma\rangle$ 
and $\langle\sigma\tau,\tau\sigma\rangle$) each act freely, and the unlabeled maps are 2-fold covering projections.
The maps $\lambda$, $g^+$, $G$ and $h$ are each 2-fold branched coverings,
since $\langle\sigma,\tau\rangle$ acts freely on $S^2\times{S^2}\setminus(\Delta_{S^2}\cup\Gamma_A)$.
The part of this diagram involving the vertices $S^2\times{S^2}$, $\mathbb{CP}^2$ ,
$\mathbb{RP}^2\times\mathbb{RP}^2$, $\mathbb{RP}^4$ and $S^4$ is displayed in \cite{Ma}.
 
On the affine piece $U_2=\{[z_0:z_1:z_2]\in\mathbb{CP}^2\mid{z_2\not=0}\}\cong\mathbb{C}^2$
we have
\[f^{-1}([p:q:1])=([\frac12(q\pm\sqrt{q^2-4p}):1],[\frac12(q\mp\sqrt{q^2-4p}):1]).\]
Hence
\[
G([p:q:1])=
\nu(4p,2p\bar{q}-2q,p\bar{p}+1-\frac12q\bar{q}-\frac12|q^2-4p)
\]
on $U_2$.
Homogenizing this formula gives
\[
G([z_0:z_1:z_2])=
\nu(4z_0\bar{z_2},2z_0\bar{z_1}-2z_1\bar{z_2},
z_0\bar{z_0}+z_2\bar{z_2}-\frac12z_1\bar{z_1}-\frac12|z_1^2-4z_0z_2|).
\]
The argument of $\nu$ is nonzero when $(z_0,z_1,z_2)\not=(0,0,0)$,
and its length is the square root of an homogeneous quartic polynomial in the 
real and imaginary parts of the harmonic coordinates of $\mathbb{CP}^2$.
The map $G$ is continuous, and is real analytic away from $\hat{f}(\Delta)=V(z_1^2-4z_0z_2)$.
Its essential structure is most easily seen after using $\theta$ to make a linear change of coordinates.
Let $\eth=G\theta$.
Then 

\medskip
\noindent
$\eth([u:v:w])=$
\[
\nu(iu\bar{u}+u\bar{w}+\bar{u}w-iw\bar{w},2(i-1)(u\bar{v}+\bar{u}v-iv\bar{w}-i\bar{v}w),
2u\bar{u}+2w\bar{w}-v\bar{v}-|2u^2+2w^2-v^2|),
\]
for all $[u:v:w]\in\mathbb{CP}^2$.
(Writing out the norm of the argument of $\nu$ explicitly 
does not lead to further enlightenment.)

It is clear from this formula that $\eth{c_2}=\eth$.  
The map $\eth$ is a 2-fold branched covering, 
with branch locus $\mathfrak{Re}(\mathbb{CP}^2)\cong\mathbb{RP}^2$,
the set of real points of $\mathbb{CP}^2$,
and so $S^4$ is the quotient of $\mathbb{CP}^2$ by complex conjugation \cite{Ku,Ma}.

We conclude by showing that the branch set $\eth(\mathbb{RP}^2)$ 
is unknotted in $S^4$.
We view $S^2$ as the unit sphere in the space $V$ of purely imaginary quaternions
and shall identify $S^3$ with the unit quaternions.
The standard inner product on $V$ is given by $v\bullet{w}=\mathfrak{Re}(v\bar{w})$,
for $v,w\in{V}$.
Let 
\[
C_x=\{(s,t)\in{S^2\times{S^2}}\mid{s\bullet{t}=x}\}, ~\forall~x\in[-1,1].
\]
Then $C_1=\Delta_{S^2}$ and $C_{-1}=\Gamma_A$,
while $C_x\cong{C_0}$ for all $|x|<1$.
The map $f:S^3\to{C_0}$ given by $f(q)=(q\mathbf{i}q^{-1},q\mathbf{j}q^{-1})$
for all $q\in{S^3}$ is a 2-fold covering projection, and so $C_0\cong{RP^3}$.
It is easily seen that $N=\cup_{x\geq0}C_x$ and $\sigma(N)$ 
are regular neighbourhoods of $\Delta_{S^2}$ and $\Gamma_A=\sigma(\Delta_{S^2})$, 
respectively.
The partition $S^2\times{S^2}=N\cup\sigma(N)$ into two pieces with common boundary $C_0$ 
is invariant under the action of $\langle\sigma^2,\tau\rangle$,
and the pieces are swapped by $\sigma$.
Hence there is an induced partition of $S^4=S^2\times{S^2}/\langle\sigma^2,\tau\rangle$ 
 into two diffeomorphic pieces.
Since the image of $\sigma(N)$ in $S^4$ is a regular neighbourhood of the branch set,
and its complement is diffeomorphic to this regular neighbourhood,
 the branch set is the image of one of the (two) standard unknotted embeddings
 of $\mathbb{RP}^2$ in $S^4$, by Theorem 3 of \cite{Pr}.

%
\medskip
{\it Remark.}
The main step in \cite{Ma} used a
result on fixed point sets of involutions of symmetric products
to obtain a diffeomorphism 
$\mathbb{RP}^2\times\mathbb{RP}^2/(x,y)\sim(y,x)\cong\mathbb{RP}^4$.
Our contribution has been the explicit branched covering 
$g:\mathbb{RP}^2\times\mathbb{RP}^2\to\mathbb{RP}^4$,
and the subsequent formulae for $G$ and $\eth$.
The argument in \cite{Ku} was very different.
Let $\eta:\mathbb{C}^3\to\mathbb{R}^6$ be the function given by
\[
\eta(z_1,z_2,z_3)=(|z_1|^2,|z_2^2|,|z_3^2|,\mathfrak{Re}(z_2\bar{z_3}),\mathfrak{Re}(z_3\bar{z_1}),
\mathfrak{Re}(z_1\bar{z_2})).
\]
Then $\eta(\zeta{v})=\eta(v)$ for all $v\in\mathbb{C}^3$ and $\zeta\in{S^1}$,
so $\eta|_{S^5}$ factors through $\mathbb{CP}^2=S^5/S^1$.
Moreover, $\eta(\bar{v})=\eta(v)$ for all $v\in\mathbb{C}^3$,
so $\eta|_{S^5}$ factors through $\mathbb{CP}^2/\langle{c_2}\rangle$.
The image of $S^5$ lies in the affine hyperplane defined by $x_1+x_2+x_3=1$.
Let $\Sigma=S^5\cap\mathfrak{Re}(\mathbb{C}^3)\cong{S^2}$.
Then $\eta|_{\Sigma}$ induces the Veronese embedding of $RP^2$ in this hyperplane.
Kuiper showed that $\eta(S^5)$ is the boundary of the convex hull of 
$\eta(\Sigma)$ in $\mathbb{R}^5$, 
and hence that $\eta$ induces a PL homeomorphism from 
$\mathbb{CP}^2/\langle{c_2}\rangle$ to $S^4$.


\end{document}